\documentclass[10pt]{article}
\usepackage{amscd,amsmath}
\usepackage{amsfonts}
\usepackage{amssymb,amsthm}
\usepackage{longtable}

\newtheorem{Def}{Definition}[section]
\newtheorem{df}[Def]{Definition}
\newtheorem{thm}[Def]{Theorem}
\newtheorem{cor}[Def]{Corollary}
\newtheorem{lem}[Def]{Lemma}

\newtheorem{rem}[Def]{Remark}

\newcommand{\im}{{\rm Im}}
\newcommand{\re}{{\rm Re}}
\newcommand{\Div}{{\rm Div}}
\newcommand{\Int}{{\rm Int}}
\begin{document}
\title{Maximal annuli with parallel planar boundaries in the 3-dimensional Lorentz-Minkowski space}
\author{\Large Juncheol Pyo}

\date{}
\maketitle
\renewcommand{\thefootnote}{}
\footnote{
Department of Mathematics,
Seoul National University,
Seoul 151-742, Korea\\
\tt{e-mail: jcpyo@snu.ac.kr}}
\smallskip

\begin{abstract}
We prove that maximal annuli in $\mathbb{L}^{3}$ bounded by circles,
straight lines or cone points in a pair of parallel spacelike planes are
part of either a Lorentzian catenoid or a Lorentzian Riemann's
example. We show that under the same boundary condition, the same
conclusion holds even when the maximal annuli have a planar end.
Moreover, we extend Shiffman's convexity result to maximal
annuli but by using Perron's method we construct a maximal annulus
with a planar end where Shiffman type result fails.
\end{abstract}

\pagestyle{plain}

\bigskip

\section{Introduction}
\bigskip
$\quad$ In 1956, Shiffman \cite{mm} proved simple but beautiful
theorems on minimal surfaces lying between two horizontal planes.
Let $M$ be a minimal annulus in $\mathbb{R}^3$, $P_1$, $P_2$
horizontal planes such that $\partial M=C_1 \cup C_2$, $C_{i}\subset
P_{i}$, $i=1,2$. First, Shiffman's circle theorem: For any
horizontal plane $P$ between $P_1$ and $P_2$, $M\cap P$ is a circle
whenever $C_{1}$, $C_2$ are circles. Second, Shiffman's convexity
theorem: For any horizontal plane $P$ between $P_1$ and $P_2$,
$M\cap P$ is a convex Jordan curve whenever $C_{1}$, $C_2$ are
convex Jordan curves.  Fang \cite{mfang2} generalized Shiffman's
circle theorem when a minimal annulus is bounded by a circle and
a straight line in parallel planes.  In case both curves are straight lines, it must be
parallel. Not only that, but Fang and Wei \cite{mfw} proved that a minimal
annulus with one planar end, bounded by straight lines or circles in
a pair of parallel planes, is part of a Riemann's example.
On the other hand, the Shiffman's convexity theorem does not hold when minimal annuli have a planar end. Fang and Hwang \cite{mfh} construct a minimal annulus with a planar end and bounded
by a circle and a strictly convex non-circular Jordan curve, but
there exists a non convex Jordan curve as an intersection curve with
a parallel plane.

By the analogy with minimal surfaces in $\mathbb{R}^3$,
F. L\'{o}pez, R. L\'{o}pez and Souam \cite{mlls} proved that (i) only
Lorentzian catenoids and Lorentzian Riemann's examples are foliated by
pieces of circles in parallel planes. And they proved a similar
result of Enneper's (see \cite{mn}), that is, (ii) a maximal
spacelike surface foliated by pieces of circles, the planes
containing these pieces of circles must be parallel. By (ii), we can
rewrite (i) as follows: Only Lorentzian catenoids and Lorentzian
Riemann's examples are foliated by pieces of circles. Also, they
proved a similar Shiffman's circle result, (iii) a maximal annulus
bounded by two parallel planar circles, then the intersection of the
annulus by a plane parallel to the boundary circle is again a
circle. Hence the maximal annulus is part of
either a Lorentzian catenoid or a Lorentzian Riemann's example.

In this paper inspired by Shiffman, Fang, Hwang and Wei's
works, we extend them to the maximal version. We have organized the
present paper as follows.

In section 2, we review some well-known
facts and previous properties on the Lorentz-Minkowski space. In
particular, we refer the F. L\'{o}pez, R. L\'{o}pez and Souam's result.

In section 3, we consider maximal annuli bounded by parallel planar
curves which are constant curvature. We prove the Lorentzian Shiffman's circle theorem.

\noindent\textbf{Theorem 1.}(Theorem \ref{eq:m6})
\textit{A maximal annulus, bounded by straight lines, circles or cone point in a
pair of parallel planes, is part of either a Lorentzian catenoid or a Lorentzian Riemann's example.
 If both curves are straight lines, they must be parallel.}

In section 4, we consider maximal annuli with a planar end and bounded by parallel planar
curves which are constant curvature.

\noindent\textbf{Theorem 2.}(Theorem \ref{eq:m11}) \textit{A maximal annulus with a
planar end, bounded by straight lines, circles or cone point in a
pair of parallel planes, is part of a Lorentzian Riemann's example.}

In section 5, we show that if a maximal annulus has a planar end then the
Lorentzian Shiffman's convexity theorem does not hold. More precisely, we
have the following theorem:

 \noindent\textbf{Theorem 3.}(Theorem \ref{eq:m14}) \textit{We
construct a maximal annulus with a planar end and it is bounded by a circle
and a strictly convex non-circular Jordan curve, but there exist a
non-convex Jordan curve as an intersection curve with a parallel plane.}
\section{Preliminaries}
$\quad$ Let $\mathbb{L}^3$ be the three dimensional
Lorentz-Minkowski space, that is, the real vector space
$\mathbb{R}^3$ endowed with the Lorentz-Minkowski metric $\langle,\rangle$,
where $\langle,\rangle = d{x_1}^2 +d{x_2}^2 - d{x_3}^2$ and $x_1, x_2, x_3$
are the canonical coordinates of $\mathbb{R}^3$. We say that a
vector $v \in \mathbb{R}^3 -\{0\}$ is spacelike, timelike or
lightlike if $|v|^2 = \langle v,v\rangle$ is positive, negative or zero,
respectively. The zero vector $0$ is spacelike by convention. A
plane in $\mathbb{L}^3$ is spacelike, timelike or lightlike if the
normal vector of the plane is
 timelike, spacelike, or lightlike, respectively. An immersed
 surface $\Sigma \subset \mathbb{L}^3$ is called spacelike if every tangent
 plane is a spacelike.  An immersed spacelike surface $\Sigma$ is
 called \textit{{maximal}} if the mean curvature is zero
 everywhere.

  Near a regular point of a maximal surface, a unit
 normal vector field can be considered as a $\it {Gauss~ map}$ $ N : \Sigma \rightarrow {\mathcal{H}}^2 =
  \{(x_1 ,x_2 ,x_3 ) \in \mathbb{L}^3 : {x_1}^2 +{x_2}^2
- {x_3}^2 = -1\}$, where ${\mathcal{H}}^2 $ the hyperbolic sphere in
$\mathbb{L}^3$ with constant intrinsic curvature identically $-1$.
Denote by $\overline{\mathbb{C}}$ the extended complex plane
$\mathbb{C}\cup \{\infty\}$. Let the stereographic projection
$\sigma$ for ${\mathcal{H}}^2$ be defined by:
$${\sigma : \overline{\mathbb{C}}- \{|z|= 1\} \rightarrow
{\mathcal{H}}^2,  ~~~~~ z \mapsto \left(\frac{2 \im (z)}{|z|^2 -1},
\frac{-2 \re (z)}{|z|^2 -1}, \frac{|z|^2 +1}{|z|^2 -1}\right)},$$ where
$\sigma(\infty)=(0, 0, 1)$, that is, $\sigma(z)$ is the intersection
of ${\mathcal{H}}^2$ and the line joining the point $(\re ({z}), \im
({z}), 0)$ and ``the north pole'' $(0, 0, 1)$ of ${\mathcal{H}}^2$.
It is well known that $\sigma$ is conformal in the natural manner.
${\mathcal{H}}^2$ has two connected components ${\mathcal{H}}^2_{+}
:= {\mathcal{H}}^2 \cap \{x_3 \geq 1 \}$ and ${\mathcal{H}}^2_{-} :=
{\mathcal{H}}^2 \cap \{x_3 \leq -1 \}$.

 Since $\Sigma$ is of zero mean curvature,
 the coordinate functions $x_1 , x_2 , x_3 $ are harmonic functions and
 hence it admits a Weierstrass
representation (see \cite{mk} for details):
\begin{thm}\label{eq:m1}(Weierstrass representation of maximal surface in
$\mathbb{L}^3$)\\
Any maximal spacelike surface in $\mathbb{L}^3$ is represented as
\begin{equation} \label{eq:m2}
 X(p) = \re \int^{p} \left(\frac{1}{2}(1+g^2)\eta, \frac{i}{2}(1-g^2)\eta,
 g\eta\right)= \re\int^{p}(\omega_{1},\omega_{2},\omega_{3}),~~~p\in D
\end{equation}
where ${D}$ is a domain in $\mathbb{C}$, and $\eta$(resp. g) is
holomorphic 1-form (resp. meromorphic function) on ${D}$ such that
$g^2\eta$ is holomorphic 1-form on ${D}$ and that $|g(\zeta)|\neq1$
for $\zeta \in {D}$. Moreover,\\ (a) The Gauss map $N$ is given by
$N(\zeta)= \sigma(g(\zeta))$. \\(b) The induced metric is given by
$ds = (|1-|g|^2||\eta|/2)$.\\(c) The Gauss curvature is given by
$K=\left[\frac{4|dg|}{|1-|g|^{2}|^{2}|\eta|}\right]^{2}$.
\end{thm}
\begin{rem}
Many properties of maximal surfaces are similar to minimal
surfaces. Contrary to the case of minimal surfaces, maximal surfaces
have naturally arising singularities due to the geometry of the
Gauss map. And since the Gauss curvatures of maximal surfaces are
always non-negative, so every maximal surface is stable.
\end{rem}

 By Calabi \cite{mc} (in general, see \cite{mcy}), every non planar complete maximal surface
has singularities. Hence many authors has studied intensively about
singularities (see \cite{misa2},\cite{misa},\cite{mkob},\cite{muy}).
Let $X:\mathcal{D}\rightarrow \mathbb{L}^{3}$ be a continuous map defined on an open
disk $\mathcal{D}$, $q$ be a interior point of $\mathcal{D}$, and suppose $X$ is
a maximal immersion on $\mathcal{D}-\{q\}$. Let $z$ be a conformal
parameter on $\mathcal{D}-\{q\}$ associated to the metric
$ds^{2}=\lambda^{2}(z)|dz|^2$ induced by $X$, where $\lambda(z)>0$
for any $z \in z(\mathcal{D}-\{q\})$. Define $q$ to be an
$\textit{isolated~singularity}$ of $X$ if for any sequence $\{q_n\}
\subset \mathcal{D}-\{q\}$ tending to $q$, the limit
$\lim_{n\rightarrow\infty}\lambda(z(q_n))$ vanishes. In this case,
we say that $X(\textit{D})$ is a maximal surface with a singularity
at $X(q)$. There are two kinds of isolated singularities called
branch points and conelike singularities.

In case $\mathcal{D}-\{q\}$ endowed with a induced complex
structure is conformally a once punctured disc, then $q$ is called
a \textit{branch~point}. This means that $\eta=0$ near $q$,
$\eta$ is a holomorphic 1-form of Weierstrass representation
and the surface cannot be embedded.

Suppose now that $\mathcal{D}-\{q\}$ is conformally to an annulus
$\{z \in \mathbb{C}: 0<r<|z|< 1 \}$.  If $X$ can be extended
continuously to $C_0=\{z \in \mathbb{C}: 0<r<|z|\leq 1 \}$, with
$X(\{|z|=1\})=X(q)$. In this case we call $q$ a
\textit{conelike~singularity}, $P_{0}=X(\{|z|=1\})=X(q)$ is called a
\textit{cone~point}, and the surface is embedded near the cone point.
At the cone point, maximal surfaces are naturally extended.

\begin{lem}\label{eq:m3}
(extension for a cone point in $\mathbb{L}^{3}$ \cite{misa})\\
Let $X_0: C=\{r<|z|<1\}\rightarrow \mathbb{L}^{3}$ be an embedded
maximal surface with cone point $P_{0}=X_{0}(\{|z|=1\})$, then the followings hold:\\
Let the Weierstrass data $(g,\eta)$ of $X_{0}$ satisfy that $g$ is
injective and $|g|=1$ on $\{|z|=1\}$ and $\eta\neq0$ on $\{|z|=1\}$.
The surface $X_{0}$ reflects analytically about
$\{|z|=1\}$ to the \textit{mirror surface}. More precisely, let
$J(z)=1/\overline{z}$ denote the inversion about $\{|z|=1\}$, the
mirror surface $X_{0}^{*}$ has the Weierstrass data
$(J^{*}g=1/\overline{g},J^{*}\phi=-\overline{\phi})$ and satisfies
$X_{0}^{*}=-X_{0}+2P_{0}$, where $P_{0}=X_{0}(\{|z|=1\})$. Moreover,
for any spacelike plane $\Pi$ contains $P_{0}$ the Lorentzian
orthogonal projection $\pi:X_{0}\rightarrow\Pi$ is a local
homeomorphism and near $P_{0}$, $X_{0}$ is asymptotic to the half
light cone with vertex at $P_{0}$.
\end{lem}
 A \textit{circle} in $\mathbb{L}^{3}$
is defined to be a planar curve with  nonzero constant curvature.
Therefore, there are three different types of circles in
$\mathbb{L}^{3}$ since there are three different types of planes in
$\mathbb{L}^{3}$. In this paper, however, circles are the same as in
$\mathbb{R}^{3}$ since we focus only on spacelike planes in
$\mathbb{L}^{3}$.  \textit{Straight lines} in $\mathbb{L}^{3}$ are
defined as similarly.

 We introduce Lorentzian Riemann's examples.
\begin{thm}\cite{mlls}\label{eq:m4}
Let $X:M\rightarrow \mathbb{L}^3$ be a spacelike conformal
non-planar maximal immersion of a Riemann surface $M$. If $X(M)$ is
foliated by pieces of Euclidean circles in parallel planes with
normal Euclidean vector $v=(0,0,1)$, then, up to scaling and linear
isometries in $\mathbb{L}^3$, the Gauss map $g$ of $X$ satisfies:
\begin{enumerate}
\item $\frac{d g}{dz}=g$, or
\item $\big(\frac{d g}{dz}\big)^2=g(g^2+2rg+1)$, where $r \in \mathbb{R}$.
\end{enumerate}
\end{thm}
We call that the first case is a \textit{Lorentzian catenoid}, the
second case is a \textit{ Lorentzian Riemann's example}.

  Now we consider a connected component of outside of Euclidean ball.
  This connected component conformally equivalent to punctured disk and the metric
has a pole at the puncture. The connected component called an
\textit{end}. The asymptotic behaviour of an end is similar to an end of minimal
surfaces (see \cite{mn} for details). Similar result in the Lorentzian
setting can be founded in \cite{mi}. Also a different approach to an
end by Klyachin can be founded in \cite{mkl}. We omit the proof.
\begin{lem}\label{eq:m5}
Let $X:\mathcal{D}\setminus\{0\}\rightarrow \mathbb{L}^{3}$ be an
embedded end of maximal surface with vertical limit normal and  the
Weierstrass data $(g,\eta)$, then the followings hold:\\
The order of pole of $\omega_{i},~i=1,2,3$ is two and
the end $X$ is asymptotic to the following:
$$(x_{1},x_{2},x_{3})=(\alpha r^{-1}\cos\theta, \alpha r^{-1}\sin\theta, \beta\log r),$$
on a neighborhood of $0$, where $z=re^{i\theta}$, $\alpha\in \Bbb R\setminus\{0\}$, and $\beta\in \Bbb R$.
\end{lem}
An end is called \textit{planar end} (resp. \textit{catenoidal end})
if $\beta=0$ (resp. $\beta\neq0$) and it is asymptotic to a horizontal plane
(resp. a vertical half Lorentzian catenoid).
\section{Maximal annuli in a slab}
By Lorentzian isometry we can denote a spacelike plane
$\Pi=\Pi_{t}=\{(x_{1},x_{2},x_{3})|x_{3}=t\}$ and a slab
$S(a,b)=\{(x_{1},x_{2},x_{3})|a\leq x_{3}\leq b\}$. By homothety
we also assume that $S(a,b)=S(-1,1)$.
\begin{thm}\label{eq:m6}
Let $A\subset S(-1,1)$ be a compact maximal annulus in a slab whose
set of singularities consists of a finite (possibly empty) set of
conelike singularities. Suppose $A(1)=A\cap \Pi_{1}$,
$A(-1)=A\cap \Pi_{-1}$ are straight lines, circles or cone
points.
\begin{enumerate}
\item If both $A(1)$ and $A(-1)$ are circles then $A(t)=A\cap
\Pi_{t}$ is a circle or cone point for $-1<t<1$. In particular, $A$
is embedded and the number of cone points is at most one.
\item If $A(1)$ or $A(-1)$ is a straight line, the other one
is a circle and $A$ is embedded, then $A(t)=A\cap \Pi_{t}$ is a
circle or a cone point for $-1<t<1$.
\item If both $A(1)$ and $A(-1)$ are straight lines and $A$ is
embedded, then $A(t)=A\cap\Pi_{t}$ is a circle or a cone point for
$-1<t<1$.
\item If $A(1)$ is a straight line or a circle and $A(-1)$ is a cone point
(in case $A(1)$ is a straight line, we also assume that $A$ is embedded) then $A(t)=A\cap
\Pi_{t}$ is a circle or a cone point for $-1<t<1$.
\end{enumerate}
\end{thm}

In order to prove the Theorem \ref{eq:m6}, we need some lemmas.

\begin{lem}\label{eq:m7}
Let $A\subset S(-1,1)$ be a properly immersed maximal annulus such that
both $A(1)$ and $A(-1)$ are circles or straight lines, then
$A$ can be conformally parameterized by
$$X:A_{R}-C\rightarrow \mathbb{L}^{3},$$
where $A_{R}=\{z \in \mathbb{C}:1/R\leq|z|\leq R\}$
for $1<R< \infty$ and the set $C$ is determined as follows:\\
If $A(1)$ and $A(-1)$ are both circles then $C=\emptyset$; if $A(1)$
is a straight line and $A(-1)$ is a circle (resp. $A(1)$ is a circle
and $A(-1)$ is a straight line) then $C=\{p : |p|=R\}$ (resp.
$C=\{q: |q|=1/R\}$); if $A(1)$ and $A(-1)$ are both straight lines
then $C=\{p,q: |p|=R, |q|=1/R\}$.\\
 In any case, the Gauss map $g$ of $A$ has neither zero nor pole in
 the interior of $A_{R}$, and it can be extended to a neighborhood of $A_{R}$. Moreover,
  the extended $g$ has either zero or pole order two at $p$ and $q$.
\end{lem}

\begin{proof}
Since $A$ is a proper maximal annulus, the conformal structure of
the interior of $A$ is equivalent to interior of $A_{R}=\{z \in
\mathbb{C}:1/R\leq|z|\leq R\}$ for some $1<R< \infty$, and a
conformal harmonic immersion $X:A_{R}-C\rightarrow S(-1,1)$, where
$C$ is a subset of $\partial A_{R}$ and $A(\{|z|=R\}-C)=A(1)$,
$A(\{|z|=1/R\}-C)=A(-1)$. In particular, the third coordinate
function $X_{3}$, which is harmonic with $X_{3}|_{(\{|z|=R\}-C)}=1$,
$X_{3}|_{(\{|z|=1/R\}-C)}=-1$ and $-1< X_{3}|_{\Int(A_{R})}<1$,
can be extended to whole $A_{R}$ such that $X_{3}|_{\{|z|=R\}}=1$,
and $X_{3}|_{\{|z|=1/R\}}=-1$. By the existence and uniqueness of
the Dirichlet problem, $X_{3}=\frac{1}{\log R}\log |z|$, we have
for any $-1<t<1$, $A(t)=A \cap \Pi_t$ is the image
$X(\{z\in A_{R}:|z|=R^t\})$.

First, $g$ cannot have zeros or poles in Int$(A_R)$, interior of $A_R$.
Suppose not, the preimage of $A(t)$ for $t$ has at least
four rays at a zero or a pole by the harmonicity of maximal surfaces.
But the preimage of $A(t)$ is a circle since $X_{3}=\frac{1}{\log
R}\log{|z|}$. So there are no zeros and poles in the interior.

It remains only to prove that on the boundary of $A$, that is, the Gauss map $N$ is not
 perpendicular to the $x_{1}x_{2}=xy$-plane. Since boundaries are composed with a circle or a straight line, the projection of the boundary into the $xy$-plane satisfies the sphere
 condition, inner or outer. There is well-defined normal direction at every
 boundary point. Near any boundary point $p$, $N$ has a vertical normal, the surface is a
 graph over a small open disk $D \subset P_1$ with $p$ on $\partial
 D$, assuming that $p\in A(1)$. Then we can write by the maximal surface
 equation. We write $(x, y, z=x_{3}) \in A$, where $x_{3}=z(x, y)$ satisfies
 $$(1-z^{2}_{y})z_{xx} +2z_x z_y z_{xy} +(1-z^{2}_{x})z_{yy}=0,~~~z_{x}^{2} +z_{y}^{2} < 1.$$
 Since $X_3$, the third coordinate function of $A$, is harmonic, by
 maximum principle we have for any $(x, y) \in D$ that $z(x,
 y)<1=z(p)$. Define a uniformly elliptic operator on a smaller
 domain if necessary,
 $$Lu=(1-u^{2}_{y})u_{xx} +2u_x u_y u_{xy} +(1-u^{2}_{x})u_{yy},~~~u_{x}^2+u_{y}^2<1.$$
 Then $z$ satisfies $Lz=0$. By Hopf boundary point lemma
 $$\frac{\partial z}{\partial\nu}>0,$$
 where $\nu$ is the outward normal to $\partial D$ at $p$. But this
 means that the normal is not vertical. This contradiction proves
 that $N$ is never vertical on the boundary of $A$. Hence $g\neq 0$
 or $\infty$.

 If $A(1)$ is a straight line, by Lorentzian isometry we can assume
 that $A(1)$ is parallel to the $y$-axis in $\mathbb{L}^{3}$. Then
 the normal vector of $A$ along the $A(1)$ stays in the $xz$-plane.
 Let $C_{1}=C\cap \{|z|=R\}$, $g$ is real on $\{|z|=R\}-C_{1}$.
 Using the Schwarz reflection principle, $g$ can be
 extended to $\{R<|z|<R^{3}\}$ by
 $\widetilde{g}(z)=\overline{g(R^{2}/\overline{z})}$ for
 $R<|z|<R^{3}$. So we get a maximal surface
 $$\mathcal{A}=\overline{X}:\{1/R<|z|<R^{3}\}-C_{1}\rightarrow S(-1,3).$$
 Since $X$ is properly immersed, the extended surface
 $\overline{X}(\{1/R<|z|<R^{3}\}-C_{1})$ is also properly immersed
 and contains a complete maximal annular end. Since the Gaussian curvature of maximal
 surface is always nonnegative, by Huber's theorem
(\cite{mh}, or see appendix of \cite{muy}) the annular end of
$\mathcal{A}$ conformally equivalent to a punctured disk
  and the Gauss map of $\mathcal{A}$ can be extended to the
  puncture. Hence $C_{1}=\{p\}$ is singletone and $g$ is either
  zero or infinite, unless the length of the straight line is finite.
  Hence $\mathcal{A}$ has a vertical limit end, by lemma \ref{eq:m5},
  at the $p$ has zero of order two.  If $A(-1)$ is
  a straight line, we apply the same process.
 \end{proof}

 Now, we derive the Lorentzian Shiffman function in terms of Weierstrass
 data. First we calculate planar curvature of each $A(t)=A\cap
 \Pi_{t}$, $-1\leq t\leq 1$. At any point of $A(t)$, let $\psi$ be
 the angle between the tangent vector and the positive $x$-axis. By
 lemma \ref{eq:m7}, $g\neq 0,\infty$ in the interior of $A_{R}$, so
 the unit normal vector is $g/|g|$, and $\phi=\arg g= \im (\log
 g)=\psi-\pi/2$. We note that the function $\phi$ can be multivalued
 but harmonic. Now suppose that $s$ is the arclength parameter of
 the curve $A(t)$, and $X^{-1}(A(t))=\{z:|z|=r=R^{t}\}$, write
 $z=re^{i\theta}$, then the curvature of $A(t)$ is:
 \begin{eqnarray} \label{eq:n1}
 \kappa&=&\psi_{s}=\phi_{s}=\frac{d}{ds}\im \left(\log
g\right)=\im\left(\frac{d}{ds}\log g\right)\\
&=&\im\left(\frac{g'}{g}\frac{dz}{d\theta}\frac{d\theta}{ds}\right)=\im
\left(\frac{g'}{g}izr^{-1}\Lambda^{-1}\right)=r^{-1}\Lambda^{-1}\re\left(z\frac{g'}{g}\right).\nonumber
\end{eqnarray}
\noindent Here we use that on the curve $\{z:|z|=r=R^{t}\}$,
$$\frac{dz}{d\theta}=ire^{i\theta},~ds=\Lambda|dz|=r\Lambda d\theta.$$
By direct calculating, we have the Lorentzian Shiffman function:
\begin{eqnarray} \label{eq:n2}
u:=r\Lambda\frac{\partial \kappa}{\partial\theta}=
\im\left[\frac{1}{2}\frac{|g|^{2}+1}{|g|^{2}-1}
\left(z\frac{g'}{g}\right)^{2}-z\frac{d}{dz}\left(z\frac{g'}{g}\right)\right].
\end{eqnarray}
\begin{lem}\label{eq:m8}
Let $A$ and $C$ be as in lemma \ref{eq:m7}, and $u$ be the
Lorentzian Shiffman function as (\ref{eq:n2}). Then $u$ can be
continuously extended on the set $C$ and $u=0$ on the boundary
$\partial A$.
\end{lem}
\begin{proof}
Let
$U(z)=\left[-\frac{1}{2}\left(z\frac{g'(z)}{g(z)}\right)^{2}-z\frac{d}{dz}\left(z\frac{g'(z)}{g(z)}\right)\right]
+\left[\left(1-\frac{1}{1-|g|^2}\right)\left(z\frac{g'(z)}{g(z)}\right)^2\right]=\Phi(z)+\Psi(z).$ We
claim that both $\Phi$ and $\Psi$ are $C^{\infty}$ complex function
near any point of the set $C$. The claim is proved then since $u(z)=\im U(z)$
is smooth near $z_{0}$, $u(z)$ can be continuously extended to $p$.

Let $z_{0}=p$ or $q$. By the lemma \ref{eq:m7}, the
extended Gauss map $\widetilde{g}$ has a zero or a pole of order two.
Let us assume that $g(z_{0})=0$.

First, we show that $\Phi$ is a $C^{\infty}$ complex function near
each point of the set $C$. Let $\zeta=z-z_{0}$, we have
$$\widetilde{g}(z)=(z-z_{0})^{2}h(z)=\zeta^{2}h(z_{0}+\zeta),$$
where $h$ is a holomorphic function and $h(z_{0})\neq0$. For
convenience, write $g$ instead of $\widetilde{g}$, then
$$z\frac{g'(z)}{g(z)}=\frac{2z_0}{z-z_0}+2+z\frac{h'(z)}{h(z)}
=\frac{a_{-1}}{\zeta}+\sum^{\infty}_{k=0}a_{k}\zeta^{k},~a_{-1}=2z_{0}.$$
And we have
$$\left(z\frac{g'(z)}{g(z)}\right)^{2}=\frac{a_{-1}^{2}}{\zeta^{2}}+\frac{2a_{-1}a_{0}}{\zeta}
+\sum^{\infty}_{k=0}b_{k}\zeta^{k},$$ and
$$z\frac{d}{dz}\left(z\frac{g'(z)}{g(z)}\right)=-\frac{a_{-1}z_{0}}{\zeta^{2}}-
\frac{a_{-1}}{\zeta}+(\zeta
+z_0)\sum^{\infty}_{k=1}ka_{k}\zeta^{k-1}.$$

Since $a_{-1}=2z_{0}$, we have $a_{-1}^{2}-2a_{-1}z_{0}=0$.
Equation $$z\frac{h'(z)}{h(z)}=(\zeta
+z_{0})\sum^{\infty}_{k=0}\frac{1}{k!}\left(\frac{h'}{h}\right)^{(k)}(z_{0})\zeta^{k}=
z_{0}\frac{h'(z_0)}{h(z_0)}+\sum^{\infty}_{k=1}c_{k}\zeta^{k},$$
implies $$a_{0}=2+z_{0}\frac{h'(z_{0})}{h(z_{0})}.$$ Now we
calculate the value $a_{0}$. The Weierstrass representation for the
extended surface $\mathcal{S}$ is
$$\phi_{1}=\frac{1}{\log R}\frac{1}{2z}\left(\frac{1}{\widetilde{g}}+\widetilde{g}\right)dz,
~\phi_{2}=\frac{1}{\log
R}\frac{i}{2z}\left(\frac{1}{\widetilde{g}}-\widetilde{g}\right)dz,~\omega_{3}=\frac{1}{\log
R}\frac{1}{z}dz.$$ For simplicity, we write $g$ instead of
$\widetilde{g}$. Let us choose a loop $\gamma$ around $z_{0}$ small
enough so the inside of $\gamma$ has only one element of set $C$. By
well-definedness of a extended maximal surface and
$$\int_{\gamma}(\omega_{1},\omega_{2},\omega_{3})=\overrightarrow{0},$$
we have $$\int_{\gamma}\frac{1}{zg}=\int_{\gamma}\frac{g}{z}=0.$$
Then
\begin{eqnarray}\label{eq:m9}
 0&=&\lim_{z\rightarrow z_{0}}\left(\frac{(z-z_{0})^2}{z g(z)}\right)'=\lim\left(\frac{1}{z h(z)}\right)'\nonumber\\
 &=&-\frac{1}{z_{0}^{2}h(z_{0})}-\frac{h'(z_{0})}{z_{0}h^{2}(z_{0})}\nonumber.
\end{eqnarray}
Finally we have
$$a_{0}=2+z_{0}\frac{h'(z_{0})}{h(z_{0})}=1.$$
Hence
$$\Phi(z)=-\frac{1}{2}\frac{a^{2}_{-1}-2a_{-1}z_{0}}{\zeta^{2}}-\frac{a_{-1}a_{0}-a_{-1}}{\zeta}
-\frac{1}{2}\sum^{\infty}_{k=0}b_{k}\zeta^{k}-(\zeta+z_{0})\sum^{\infty}_{k=0}ka_{k}\zeta^{k-1}\\
=\sum^{\infty}_{k=0}d_{k}\zeta^{k}$$ is holomorphic near
$z=z_{0}$.

Now we consider the function $\Psi(z)$. Since
$|g(z)|^{2}=|z-z_{0}|^{4}|h(z)|^{2}=|\zeta|^{4}|h(z)|^{2}$ and
$\zeta^{2}\left(z\frac{g'(z)}{g(z)}\right)^{2}$ is holomorphic, it follows that
$$\frac{1}{\zeta^{2}}\left(1-\frac{1}{|g|^{2}-1}\right)=\frac{1}{\zeta^{2}}\sum^{\infty}_{k=1}|g|^{2k}
=\zeta^{2}\sum^{\infty}_{k=1}|\zeta|^{4(k-1)}|h(z)|^{2k}$$ is a
smooth function near the $z_{0}$. Thus $\Psi(z)$ is a smooth
function near $z_{0}$, and so $U(z)$ is also smooth.
Since $u|_{ \partial A-C}=0$ and $u$ can be continuously extend to
$C$, $u=0$ on the $\partial A$.
\end{proof}
\begin{lem}\label{eq:m10}
The Lorentzian Shiffman function $u$ can be smoothly extended on the
 conelike singularities.
\end{lem}
\begin{proof}
Let $X:\{r<|z|<1\}\rightarrow \mathbb{L}^3$ has cone point at
$X(\{c:=|z|=1\})$ with Weiersrtrass data $(g,\phi)$. By M\"{o}bius
transformation on $c$, we can assume that the curve $c$ is
$Re(z)=0$. And the involution $J$ is $J(z)=-\overline{z}$, and the
Weierstrass data of mirror surface are
$(1/\overline{g},-\overline{\phi})$. Write
$g=e^{w(z)}$, we have $w(-\overline{z})=-\overline{w(z)}$.\\
The Lorentzian Shiffman function $u$ extends to $c$ if and only if
$$V(z):=\im\left(\frac{1}{|g^2-1|}\left(\frac{d\log g}{dz}\right)^2\right)$$ extends to
$c$.

we claim that $V(z)$ can be smoothly extended to $c$.

\noindent Take $z_{0}\in i\mathbb{R}$ and let the Taylor series of the function
$w(z)=\sum^{\infty}_{m=0}a_{m}(z_{0})(z-z_{0})^{m}$. Since
$w(-\overline{z})=-\overline{w(z)}$, we have
$(-1)^{m}a_{m}(z_{0})=-\overline{a_{m}(z_{0})}$, that is,
$\re(a_{2n}(z_{0}))=0$ and $\im(a_{2n+1}(z_{0}))=0$ for all $n \in
\mathbb{N}\cup\{0\}$. Since $g$ is injective near a conelike
singularity, we
have $a_{1}(z_{0})\neq0$, for any $z_{0}\in c$.

Now we have
$$|g(z)|-1=e^{\re(w(z))}-1=\re(w(z))\widetilde{H}_{1}(\re(w(z))),$$
where $H_{1}(z)=(e^{z}-1)/z$, $z\in \mathbb{C}$. Since the
coefficient of function $w(z)$, we deduce $\re w(z)=\re
(z)V_{1}(z)$, where $V_{1}$ is a suitable smooth function around
$c$, and because of $a_{0}(z_{0})\neq0 $ for all $z_{0}$ in compact
set $c$, $|V_{1}|_{c}\geq \epsilon>0$. Thus, $$|g(z)|-1=\re
(z)H_{1}(z),$$ where smooth function $H_{1}(z)$ with
$|H_{1}|_{c}\geq \epsilon'>0$. By similar argument, we have
$$\im\left(\left(\frac{d\log g}{dz}\right)^2\right)=\re (z)H_{2}(z),$$ where $H_{2}$ is
a smooth function around $c$. Hence around $c$,
$$V(z)=\frac{H_{2}(z)}{H_{1}(z)(1+|g(z)|)}$$ is a smooth function.
\end{proof}

\textit{Proof of Theorem \ref{eq:m6}.} Let us show that
$A:M\rightarrow \mathbb{L}^{3}$ bounded by two cone points $P_{1}$
and $P_{-1}$ is not possible. If not by lemma 2.3 successive
reflections about cone points, we have a complete maximal annulus
$\widetilde{A}:\widetilde{M}\rightarrow \mathbb{L}^{3}$ with
infinitely many conelike singularities such that $\widetilde{A}$ is
a translation invariant. The quotient of $\widetilde{M}$ under the
holomorphic translation induced by above translation gives a torus
T, and the Weierstrass data $(\omega_{1},\omega_{2},\omega_{3})$ of
$\widetilde{A}$ can be induced on $T$. Furthermore, $\omega_{j}$ is
holomorphic, and so $\omega_{j}=\lambda_{j}\tau_0$, $j=1,2,3$, where
$\lambda_{j}\in \mathbb{C}$ and $\tau_0 $ is a nonzero holomorphic
1-form on $T.$ Because
$\omega_{1}^{2}+\omega_{2}^{2}-\omega_{3}^{2}=0$ and the associated
maximal immersion is singly periodic, it is not hard to see that
$\lambda_j=r_j\lambda,$ where $r_j \in \mathbb{R}$, $\lambda \in
\mathbb{C}$ and $r_{1}^{2}+r_{2}^{2}-r_{3}^{2}=0$. In particular,
$\widetilde{A}$ lies in a lightlike straight line, which is
impossible.

First, both $A(1)$ and $A(-1)$ are circles. By lemma \ref{eq:m7}, we
find a conformal annulus $A_R$ and the set $C$ is empty. By lemma
\ref{eq:m8} and \ref{eq:m10}, the Lorentzian Shiffman function $u$
is a smooth in the interior of $A_R$ and $u$ satisfies
$$\triangle_{A} u=2K u,~~u|_{\partial A_{R}}=0.$$
Since every maximal surface is stable, the first eigenvalue of
Jacobi operator is positive. Hence $u\equiv0$ and $\Pi_{t}$ is a
circle or a conelike singularity, for $-1<t<1$. Moreover, $A$ is
part of the Lorentzian catenoid or a Lorentzian Riemann's example. So the maximal
annulus is embedded. Because Lorentzian Riemann's examples can have
at most one cone point without planar end, the maximal annulus has
at most one cone point.

Second, $A(1)$ is a straight line and $A(-1)$ is a circle. By lemma
\ref{eq:m7}, the function $u$ smooth near $A(1)$ and zero on the
$A(1)$. The same argument in the first case still holds. The third case is similar to the
second case.

Finally, either $A(1)$ or $A(-1)$ is cone point and the other is
circle or straight line. Using the lemma \ref{eq:m3}, we obtain
maximal annulus bounded by circles or straight lines. So it is a
previous case. The theorem is complete.  $\Box$

\begin{cor}
Let $A(1)$ and $A(-1)$ be non parallel straight lines to each
other. Then $\Gamma=A(1)\cup A(-1)$ cannot bound a properly embedded
maximal annulus in $S(-1,1)$.
\end{cor}
\section{Maximal annuli with a planar end in a slab I}
In this section, we consider maximal annuli with an end. This
gives a characterization of Lorentzian Riemann's examples.
\begin{thm}\label{eq:m11}
Let $A\subset S(-1,1)$ be an embedded maximal annulus with a
planar end in a slab whose set of singularities
consists of a finite (possibly empty) set of conelike singularities.
Suppose $A(1)=A\cap \Pi_{1}$, $A(-1)=A\cap \Pi_{-1}$ are
straight lines, circles or cone points, except bounded by two cone
points then $A(t)=A\cap \Pi_{t}$ is a circle or cone point for
$-1<t<1$, except at the height of the end where the intersection is
a straight line. Consequently, $A$ is part of a Lorentzian Riemann's
example, so if the boundary consists of two straight lines then the
lines must be parallel.
\end{thm}

\begin{lem}\label{eq:m12}
Let $A\subset S(-1,1)$ be a maximal annulus with a planar end and
both $A(1)$ and $A(-1)$ consist of circles or straight lines, then
$A$ can be conformally parameterized by
$$X:A_{R}-C-\{z_{e}\}\rightarrow \mathbb{L}^{3},$$
where $A_{R}=\{z \in \mathbb{C}:1/R\leq|z|\leq R\}$
for $1<R< \infty$ and the set $C$ determined as follows:

For $|p|=R$ and  $|q|=1/R$, $C=\{p,q\}$ if $A(1)$ and $A(-1)$ are
straight lines; $C=\{p\}$ (resp. $C=\{q\}$) if $A(1)$ is a straight
line and $A(-1)$ is a circle (resp. $A(1)$ is a circle
and $A(-1)$ is a straight line); and $C=\emptyset$ otherwise.

In any case, the Gauss map $g$ of $A$ has neither zero nor pole in
the interior of $A_{R}$, and $g$ can be extended to a neighborhood
of $A_{R}$ such that the extended $g$ has either zero or pole of
order two at $z_{e}$, $p$ and $q$.
\end{lem}
\begin{proof}
Since the Gaussian curvature of a maximal surface is always
nonnegative. By Huber's theorem, the conformal domain of a maximal
surface is $A_R-C-\{z_e\}$. By lemma \ref {eq:m5}, the gauss map $g$
has zero or pole of order two at $z_e$. For the rest parts of lemma are
proved the same way as in lemma \ref {eq:m7}.
\end{proof}
\textit{Proof of Theorem \ref{eq:m11}.} Either $A(1)$ or $A(-1)$ is
a cone point, using lemma \ref{eq:m3}, the maximal annulus can be
extended to maximal surface bounded by circles or straight lines. By
the lemma \ref{eq:m8}, lemma \ref{eq:m10} and lemma \ref{eq:m12},
the Lorentzian Shiffman function $u$ can be smoothly extended to the
set $C$, the end $z_e$ and cone points. And $u$ satisfies $$\triangle_{A}
u=2K u,~~u|_{\partial A_{R}}=0.$$ By the stability of Jacobi operator, $u\equiv 0$. So
the theorem is complete. $\Box$
\section{Maximal annuli with a planar end in a slab II}
First, we extend the Shiffman's convexity theorem to Lorentzian
space.
\begin{thm}\label{eq:m13}
Let $A\subset S(-1,1)$ be a properly immersed maximal annulus where
$A(1)$ and $A(-1)$ consist of convex Jordan curve then $A\cap
\Pi_{t}$ is a strictly convex Jordan curve for every $-1<t<1$. In
particular, $A$ is embedded.
\end{thm}
\begin{proof}
Let the angle function $\psi$, the planar curvature $\kappa$ as
(\ref{eq:n1}). Define $h=\re\left(z\frac{g'}{g}\right)=r\Lambda\kappa$,
$r\Lambda>0$. the $h$ is a harmonic and non-negative on the
boundary. By the strong maximum principle, the $h$ is strictly positive.
Thus $A\cap \Pi_{t}$ is locally strictly convex. Similar to a
minimal surface \cite{mfang}, the period of the angle function
$\psi$ is exactly $2\pi$. Hence $A\cap \Pi_{t}$ is strictly
convex.
\end{proof}
\begin{thm}\label{eq:m14}
We construct a maximal annulus $A\subset S(-1,1)$ with a
planar end and it satisfies the following properties: For some $t_{0}\in
(-1,1)$, $A\cap \Pi_{t_{0}}$ is a non-convex Jordan curve even when
the boundary $\partial A$ consists of a circle and a strictly convex
real analytic Jordan curve.
\end{thm}
\begin{lem}\label{eq:m15}
Let $A\subset S(-1,1)$ be a maximal annulus with a planar end, and
the boundary $\partial A$ consists of two Jordan curves lying a pair of
parallel planes which are boundary of $S(-1,1)$, then $A$ can be conformally parameterized by
$$X:A_{R}-\{z_{e}\}\rightarrow \mathbb{L}^{3},$$ where $A_{R}=\{z \in \mathbb{C}:1/R\leq|z|\leq R\}$
for $1<R< \infty$ and $1/R<|z_{e}|< R$.

Moreover, if $A\cap \Pi_{t}$ are strictly convex $C^{2}$ Jordan
curves for $-1\leq t \leq1$ except at $t_{0}\in(-1,1)$, the height
of the end, then $A \cap \Pi_{t_{0}}$ is a straight line.
\end{lem}
\begin{proof}
Because $A$ is a maximal annulus with a planar end, in the interior,
as lemma \ref{eq:m12} $A$ is conformally equivalent to
$A_{R}-\{z_{e}\}$ for suitable $1<R< \infty$ and $1/R<|z_{e}|< R$.
And the same argument Dirichlet problem $X_{3}\equiv\frac{1}{\log
R}\log |z|$, the planar curvature
$\kappa(z)=|z|^{-1}\Lambda^{-1}\re(z\frac{g'}{g})$. Since
$z_{e}\neq0$, ${g'}/{g}$ is  meromorphic and has an isolated
pole at $z_{e}\neq0$, $\Xi(z)=\re(z\frac{g'}{g})=|z|\Lambda\kappa(z)$ takes
positive and negative values near $z_{e}$. So $\Xi^{-1}(0)$ is a non
empty set and a real analytic 1-dimensional variety except isolated
points $\{p_{i}\}\subset \Xi^{-1}(0)$, at the $p_{i}$, $D\Xi(p_{i})$
is zero and at least four equal angular curves emit from $p_{i}$.
But $A\cap \Pi_{t}$ are strictly convex except $t_{0}$, $\Xi\neq0$
except $|z|=|z_{0}|$. So $\Xi^{-1}(0)\subset \{z:|z|=t_{0}\}$,
$\Xi^{-1}(0)$ has no singularities and is a 1-dimensional manifold
without boundary. This means $\Xi^{-1}(0)=
\{z:|z|=t_{0}\}-\{z_{e}\}$. Hence
$\Xi(z)=\re(z\frac{g'}{g})=|z|\Lambda\kappa(z)\equiv0$ on
$\{z:|z|=t_{0}\}-\{z_{e}\}$. The only case is $\kappa\equiv0$, so
$X(\{z:|z|=t_{0}\}-\{z_{e}\})$ is a straight line.
\end{proof}
We are going to construct a maximal surface by solving the exterior
Dirichlet problem for maximal surface equation. The variational
problem of area functional leads to the following divergence form of
maximal surface equation:$$Q\nu=\Div
\left(\frac{D\nu}{\sqrt{1-|D\nu|^{2}}}\right),$$ with $|D\nu|\leqslant1$.
Because we will use the Perron's method, we define the subsolution
and supersolution to the maximal surface equation.
\begin{df}(subsolution and supersolution)\\
Let $\Omega$ be a domain in $(x_{1},x_{2})$ plane. A $C^{0}(\Omega)$
function $\alpha$ is a \textit{subsolution}(resp.
\textit{supersolution}) in $\Omega$ if for every ball $B\Subset
\Omega$ where $\overline{B}\subset \Omega$, and every function $\nu$
satisfying $Q\nu=0$ in $B$ and $\alpha\leqslant\nu$(resp.
$\alpha\geqslant\nu$) on $\partial B$, then we have
$\alpha\leqslant\nu$(resp. $\alpha\geqslant\nu$).
\end{df}
We follow the classical Perron's strategy (see \cite{mgt}). (i) A
subsolution (supersolution) in a domain $\Omega$ satisfies the
strong maximum principle. (ii) Let $\nu$ be a subsolution in
$\Omega$ and $B$ be a ball strictly contained in $\Omega$. Denote by
$\overline{\nu}$ the solution in $B$ satisfying $\overline{\nu}=\nu$
on the boundary $\partial B$. We define the \textit{solution
lifting} of $\nu$ in $B\Subset \Omega$ ($\overline{B}\subset \Int{\Omega}$) by
\begin{displaymath}
V(x)= \left\{ \begin{array}{ll} \overline{\nu}(x), & \textrm{$x\in
B$}\\ \nu(x) & \textrm{$x\in \Omega- B$}.\end{array}\right.
\end{displaymath}
Then the function $V$ is also subsolution in $\Omega$.
  (iii) If $\nu_{1}$ and
$\nu_{2}$ are subsolutions (resp. supersolutions) to the maximal
surface equation, using the maximum principle, $\sup\{\nu_{1},
\nu_{2}\}$ (resp. $\inf\{\nu_{1}, \nu_{2}\}$) is a subsolution
(resp. supersolution) to the maximal surface equation.

For a continuous function $\varphi$ defined on $\partial\Omega$,
define $S_{\varphi}$ to be the set of subsolutions to the maximal
surface equation which are $C^{0}(\overline{\Omega})$ and equal to
$\varphi$ on the boundary $\partial\Omega$. The guarantee of
$S_{\varphi}$ is non empty and existence of a supersolution
$\nu^{+}$, the function $\mu(x)= \sup_{\nu\in S_{\varphi}}\nu(x)$
solves the Dirichlet problem
\begin{displaymath}
 \left\{ \begin{array}{ll} Q\mu=0, & \textrm{in $
\Omega$}\\ \mu=\varphi & \textrm{on
$\partial\Omega$}.\end{array}\right.
\end{displaymath}
This is a classical argument by iterating the solution lifting in
small balls.

\textit{Proof of Theorem \ref {eq:m14}}. Choose a $C$, a Lorentzian
catenoid with cone point in $\Pi_{0}=x_{1}x_{2}$-plane.  Let
$C^{+}=C\cap \{(x_{1},x_{2},x_{3})\in \mathbb{L}^{3}:x_{3}\geq1\}$
and $D_{1}\subset \Pi_{1}$ be a disk with $\partial D_{1}=C\cap \Pi_{1}$,
$C^{-}$ be a reflection of $C^{+}$ with respect to $\Pi_{1}$.  $C^{-}$ is a
graph of the function $v:\Pi_{0}-D\rightarrow \mathbb{R}$, where $D$
is a vertical translation of $D$ to the plane $\Pi_{0}$.

Kim and Yang \cite{mky} construct maximal surfaces asymtotic to the
Lorentzian catenoid with genus $k$. Let $\mathcal{M}_{1}$ be a Kim
and Yang's example with gunus $1$. We cut $\mathcal{M}_{1}$ at a
sufficient large height. Then we gain an annular end $E$ such that
$E$ is a graph and has non circular real analytic strictly convex
boundary lying in a horozontal plane. Because $C^{-}\cap \Pi_{-1}$
is a circle with enclose a disk $D_{-1}$, we can translate the end
$E$ in a way that $\partial E\subset \Pi_{-1}\setminus D_{-1}$ and
$E\cap C^{-}=\emptyset$. Denote $B_{-1}\subset \Pi_{-1}\setminus
D_{-1}$ the closed bounded convex domain bounded by $\partial E$ in
$\Pi_{-1}$. Let $B$ be the vertical translation of $B_{-1}$ to
$\Pi_{0}$ and denote $\Omega=\Pi_{0}-(D\cup B)$. The annular end $E$
is a graph as $w:\Pi_{0}-B\rightarrow \mathbb{R}$ such that
$w\equiv-1$ on $\partial B$.

On $\Omega$, we have a subsolution (resp. supersolution)
$\nu^{-}=\sup\{v,-1\}$ (resp. $\nu^{+}=\inf\{1,w\}$) to the maximal
surface equation. They satisfy the boundary condition
$\nu^{\pm}=\varphi$ on $\partial \Omega$, where $\varphi:\partial
\Omega\rightarrow \mathbb{R}$ is the function
\begin{displaymath}
 \left\{ \begin{array}{ll} \varphi=1, & \textrm{on $
\partial D$}\\ \varphi=-1 & \textrm{on
$\partial B$}.\end{array}\right.
\end{displaymath}
Hence $\mu(x)= \sup_{\nu\in S_{\varphi}}\nu(x)$ solves the Dirichlet
problem
\begin{displaymath}
 \left\{ \begin{array}{ll} Q\mu=0, & \textrm{in $
\Omega$}\\ \mu=\varphi & \textrm{on
$\partial\Omega$}.\end{array}\right.
\end{displaymath}
So the graph of $\mu$ is a maximal surface $A$ bounded by a circle
and a non circular real analytic strictly convex Jordan curve in
$\Pi_{1}$ and $\Pi_{-1}$ respectively. Since $-1\leq\mu\leq1$,
$A\subset S(-1,1)$ with an end. So the end must be  planar. Let
$t_{0}\in (-1,1)$ be the height of the end, the intersection curve
$A\cap \Pi_{t_{0}}$ is not a straight line. Suppose not $A\cap
\Pi_{t_{0}}$ is straight line denote $A_{s}$ be a subannulus of $A$
bounded by a circle and a straight line. By the theorem \ref
{eq:m6}, $A_{s}$ is a part of Lorentzian Riemann's example, thus $A$
is also a part of Lorentzian Riemann's example. This contradict to
the boundary condition of $A$. Hence by the lemma \ref {eq:m15},
there exist $t\in (-1,1)$ such that $A\cap \Pi_{t}$ is a non convex
Jordan curve. $\Box$

\end{document}